\documentclass{article}

\usepackage{arxiv}

\usepackage[utf8]{inputenc} 
\usepackage[T1]{fontenc}    
\usepackage{hyperref}       
\usepackage{url}            
\usepackage{booktabs}       
\usepackage{amsfonts}       
\usepackage{nicefrac}       
\usepackage{microtype}      
\usepackage{lipsum}
\usepackage{cancel}
\usepackage{graphicx}
\graphicspath{ {./images/} }

\usepackage{hyperref}
\usepackage{amsmath}
\usepackage{amssymb}
\usepackage{xcolor}
\usepackage{bm}
\usepackage{parskip} 
\usepackage{optidef}
\usepackage{multirow} 
\usepackage{bbold}
\usepackage{amsthm}  

\newtheorem{theorem}{Theorem}

\usepackage{graphicx}
\usepackage{subcaption}
\usepackage{caption}

\newcommand{\MI}{\text{MI}}

\definecolor{main}{HTML}{5989cf}    
\definecolor{sub}{HTML}{cde4ff}     

\usepackage[many]{tcolorbox}    	

\usepackage{setspace}               
\usepackage{multicol}               

\newtcolorbox{boxC}{
    colback = sub, 
    boxrule = 0pt  
}

\title{Mutual Information second order expansion is the Pearson's chi-square statistic}

\author{
 Marius Marinescu \\
  Engineering School of Fuenlabrada\\
  King Juan Carlos University \\
  Madrid, Spain \\
  \texttt{marius.marinescu@urjc.es} \\
}

\begin{document}
\maketitle
\begin{abstract}

We show that MI connects subtly and elegantly the two best-known state-of-the-art independence test statistics: the $G^2$ and the Pearson's chi-square statistic $\chi^2$. Furthermore, we show that the MI connects directly those statistics by an elegant formula arising from a stochastic Taylor expansion of MI ($\delta$-method). MI second order term is precisely $\chi^2$ up to a scale factor. As a consequence, by this connection, the difference between $G^2$ and $\chi^2$ can be explicitly quantified.
\end{abstract}

\section{Introduction}

Mutual Information (MI) is a key concept arising from information theory \cite{shannon1948, thomas2006}. One of its main application is to measure dependence between random variables (a stronger concept than correlation). In addition, multidimensional version exists \cite{MImulti, watanabe1960information, wu2009MItest}. 
Although mutual information is zero if and only if two random variables are statistically independent, it is primarily regarded as a measure of statistical dependence rather than a hypothesis-testing statistic. In practice, independence is usually assessed using dedicated statistical tests, largely because reliable estimation of mutual information from finite samples is challenging. Nevertheless, several MI-based independence tests have been proposed in the literature \cite{wu2009MItest, ai2024MItest, berrett2019MItest}, particularly in the past two decades.
Similarly, 
during the past few decades, a variety of interesting methods for testing statistical independence have been developed, including tests based on energy statistics, distance covariance, reproducing kernel Hilbert spaces (HSIC), and bootstrap resampling techniques \cite{gabor2007, gabor2009, gabor2013energy, arthur2005, lopez2013}.



In this paper, we show that the canonical measure of dependence, MI, naturally connects the two best-known independence statistics: the $G^2$-test statistic (aka likelihood-ratio statistic) and Pearson's chi-square statistic. 
Specifically, we derive an explicit relationship between these statistics by applying a stochastic Taylor expansion ($\delta$ method) to MI, revealing a direct and elegant connection between information-theoretic and classical hypothesis-testing frameworks.

%
Furthermore, we unveil the exact relation between $G^2$-test statistic  and Pearson's chi-square statistic.
Although it is well known that these two statistics are asymptotically equivalent, their connection is usually described only through the leading-order approximation \cite{agresti2012}. Here, by performing a stochastic second-order Taylor expansion of MI, we express the G-test statistic as Pearson's chi-square statistic plus explicit third-order correction terms.
%


\section{Second order expansion}

Consider that we have an independent, identically distributed (i.i.d.) sample of $n$ observations from a pair of random variables $(X,Y)$. Suppose that the support is $\Omega_{X,Y}=\{(x_1,y_1), (x_2,y_1), \dots, (x_{I},y_1),(x_1,y_2), (x_2,y_2), \dots, (x_I,y_J) \}$. We will use the following short-hand notation:  $p_{ij}=P(X=x_i, Y=y_j)$, 
 $p_{i*}=\sum_{j=1}^J p_{i j}$ and $p_{* j}=\sum_{i=1}^I p_{i j}$, for $i =1,\dots ,I$, and $j =1,\dots,J$.  
 %
%
%
Since we only have a sample, and we don't know any $p_{ij}$, we have to work with relative frequencies, also known as the empirical distribution or empirical processes \cite[Chapter 19]{vaart2000}, \cite[Section 3.8]{vaart1996}. Let $n_{ij}$ denote the observed frequency of the value $(x_i,y_j)$. Then, the non-parametric estimator of the probability $p_{ij}$ is: $\hat{p}_{ij}=\frac{n_{ij}}{n}$. Let $\bm{p}$ and $\hat{\bm{p}}_n$ design the column vectorisation of the true distribution and the relative frequency, respectively.
Finally, we denote $\bm{p}^*$ the column vectorisation of the distribution constructed as the product of the marginals (independent case), and $\hat{\bm{p}}_n^*$ the column vectorisation of its maximum likelihood estimator through the sample. 

The MI is defined as:
\begin{equation}\label{eq: MI}
    \MI(\bm{p})= \sum_{i,j=1}^{I,J} p_{ij}\ln \frac{p_{ij}}{p_{i*}p_{*j}}. 
\end{equation}
The MI can be seen as a KL divergence measure between the joint probability distribution ($\bm{p}$) and the ones constructed as the (outer) product of the marginals ($\bm{p}^*$) \cite{thomas2006}. 
We write the following theorem adapted from our previous work \cite{marinescu2025, marinescu2026bias}.

\vspace{5mm}

\begin{theorem}
    Let ${\bm{p}}$ be the joint pmf of two discrete random variables and denote as $\hat{\bm{p}}_n$ its estimator through the sample. 
Let $H= \{\frac{\partial^2\MI}{\partial p_{st}p_{ij}} \}_{\substack{s,i=1,...,I;\\ t,j=1,...,J}}, {\scriptstyle  (s,t),(i,j) \neq  (I,J)},$ be the Hessian matrix evaluated at $\bm{p}$.
Finally, let $\Sigma$ be the covariance matrix of $\sqrt{n}{\hat{\bm{p}}}_n$,\  $\bm{\lambda} = \text{eig}(\Sigma H)=\text{eig}(H \Sigma)$, and let $\bm{\chi^2}$ denote a column vector of $I \cdot J-1$ i.i.d. chi-square random variables with one degree of freedom. 

If $\MI(\bm{p}) = 0$, the following result holds 

\begin{equation}\label{eq:MI2order2}
     2{n}\MI(\hat{\bm{p}}_n) =  B_n + o_\mathbb{P}(1) \ \overset{asym}{\sim} \ \bm{\lambda}^\top \begin{pmatrix}
         \bm{\chi}^2
     \end{pmatrix} 
\end{equation}
where $  B_n= \sqrt{n}(\hat{\bm{p}}_n-  \bm{p})^\top H \sqrt{n}(\hat{\bm{p}}_n-  \bm{p})  $, and $ B_n\xrightarrow{d}B  $ with  $B \sim \bm{\lambda}^\top 
\begin{pmatrix} 
    \bm{\chi}^2
\end{pmatrix} $.
Furthermore, $\mathbb{E}[B_n] = \mathrm{tr}(H\Sigma) +  o(1)$ and $\mathbb{V}[B_n] = 2\mathrm{tr}(H\Sigma H\Sigma)  + o(n^{-1})$.
\end{theorem}

See \cite{marinescu2026bias} for a demonstration and a general second order expansion for non-vanishing MI gradient.
Rewriting from \cite{marinescu2026bias}, the Hessian matrix is:

\begin{align}\label{eq:H}
    H &= A + B + C \notag \\
      &= A_{IJ-1 \times IJ-1} + \begin{bmatrix}
        1 \\
        1 \\
        \vdots \\
        1
    \end{bmatrix}
    \begin{bmatrix}
      \bm{b}_1 & \cdots & \bm{b}_1 & \bm{b}_2 
    \end{bmatrix}
    + \begin{bmatrix}
        \bm{c}_1 \\
        \vdots \\
        \bm{c}_1  \\
        \bm{c}_2
    \end{bmatrix} \begin{bmatrix}
        1 & 1 & \cdots &1
    \end{bmatrix}
\end{align}

with
\[
\begin{array}{c c}
     A_{ij, st}=\begin{cases}
        \frac{1}{p_{ij}}-\frac{1}{p_{i*}} -\frac{1}{p_{*j}} & \ \ s = i,\ t=j \\
          -\frac{1}{p_{i*}} & \ \ s = i,\ t \neq j  \\
           -\frac{1}{p_{*j}}  & \ \ s \neq i,\ t=j  \\
            0 & \ \ s \neq i,\ t \neq j
          \end{cases} &
\begin{array}{c}
    \bm{b}_1 =  \begin{pmatrix}
        \alpha & \cdots & \alpha &\beta 
    \end{pmatrix}_{I \times 1}, \bm{b}_2 =\begin{pmatrix}
        \gamma & \cdots & \gamma
    \end{pmatrix}_{I-1 \times 1} \\
    \ \ \ \ \ \ \ \ \ \ \bm{c}_1 = \begin{pmatrix}
        0 & \cdots &  0 & \frac{1}{p_{I\cdot}} 
    \end{pmatrix}_{I \times 1}^\top \hspace{-0.5mm},   \ \bm{c}_2 = \begin{pmatrix}
         \frac{1}{p_{\cdot J}}  & \cdots &  \frac{1}{p_{\cdot J}}
    \end{pmatrix}_{I-1 \times 1}^\top
\end{array}
\end{array}
\]



\begin{equation*}
    \alpha= \left( \frac{1}{p_{IJ}}-\frac{1}{p_{I*}} -\frac{1}{p_{*J}}\right), \ \ \beta=  \left( \frac{1}{p_{IJ}} -\frac{1}{p_{*J}}\right), \ \ \gamma= \left( \frac{1}{p_{IJ}}-\frac{1}{p_{I*}}\right) 
\end{equation*}

See \cite[Fig. 1]{marinescu2025} for a visual representation of the H matrix for specific distributions, which unveils a interesting matrix structure.
To compute probabilities involving $B \sim \bm{\lambda}^\top 
\begin{pmatrix} 
    \bm{\chi}^2
\end{pmatrix} $
see \cite{moschopoulos1984}. More generally, for quadratic forms of normal random variables, see \cite{Abhranil_Das2025}.

\section{Independence Hypothesis testing}\label{sec:independence}

Testing the hypothesis that a pair of random variables $X,Y$ are independent gives the null hypothesis:
\begin{equation}
H_0: p_{i j}=p_{i *} p_{* j}, \quad i=1, \ldots, I, \quad j=1, \ldots, J
\end{equation}
against the negation of it. We take advantage of the property that MI is zero if and only if the random variables are independent. Using the asymptotic distribution for the MI estimator, we may construct an independence hypothesis test based on it.
Specifically, we have that 

\begin{equation}\label{eq:MI2}
    2n(\MI(\hat{\bm{p}}_n) -0) =  n(\hat{\bm{p}}_n- \bm{p})^\top H \bigg|_{\bm{p}} (\hat{\bm{p}}_n- \bm{p}) +  o_p(1).
\end{equation}

where (eq. \ref{eq:MI2order2}) 

\begin{equation}\label{eq:chi}
     \sqrt{n}(\hat{\bm{p}}_n- \bm{p})^\top H \bigg|_{\bm{p}} \sqrt{n}(\hat{\bm{p}}_n- \bm{p})  \xrightarrow[n \to \infty]{L} \bm{\lambda}^\top \begin{pmatrix}
         \bm{\chi}^2
     \end{pmatrix}.
\end{equation}
Under independence the eigenvalues of $\Sigma H$ collapses to $0, 1$, and $\mathrm{tr}(H\Sigma) = (I-1)(J-1)$.
In that case, $ G_n^2~\triangleq~2n\MI(\hat{\bm{p}}_n)$ is asymptotically distributed as a $\chi^2_{(I-1)(J-1)}$ random variable. 
Let $ \bm{p}_{H_0} \triangleq  \text{vec}(\bm{p}_X \cdot \bm{p}_Y^\top)$, where $\bm{p}_X$, $\bm{p}_Y$ represent the theoretical marginal distribution of $X$ and $Y$, respectively. 
%
In order to perform the hypothesis test, two candidates arise:
\begin{align}
    T^1_n &= 2n(\MI(\hat{\bm{p}}_n) - \MI(\bm{p}_{H_0}))   \stackrel{H_0}{=}
 2n\MI(\hat{\bm{p}}_n)  \text{ and } \\
    T^2_n &= n(\hat{\bm{p}}_n-  \bm{p}_{H_0})^\top H \bigg|_{\bm{p}_{H_0}} (\hat{\bm{p}}_n-  \bm{p}_{H_0}).
\end{align}
See \cite[Fig. 2]{marinescu2025} for a simulation of both statistic and their sample distribution which are compared to the asymptotic one.
Notice that the difference between them is of order $o_p(1)$. In particular the difference is of the form $ 2n \cdot  R(|| \hat{\bm{p}}_n - \bm{p}||^2)$, where the remainder $R$ is composed by third order derivatives of MI, which on view of the form of H, they are (multiplicative) inverse of powers of terms $p_{ij} \ \forall i,j$.


It is interesting to observe that $T_n^1$ does not depend on the unknown values $\bm{p}_{H_0}$, nevertheless, $T_n^2$ does. 
%
%
Thus, for $T_n^1$ we may construct an hypothesis test by considering the acceptance region $ R_0 = \{ T_n^1 \in \mathbb{R}^+: T_n^1 < \chi_{\bm{\lambda}, \alpha} \}$ with $P (\chi_{\bm{\lambda}} < \chi_{\bm{\lambda}, \alpha})=1-\alpha$.
For $T_n^2$ we have that the statistic depends on the true marginal distribution $\bm{p}_{H_0}$, which is typically unknown. A practical solution is to consider a plug-in estimator and substitute  $\bm{p}_{H_0}$ by the product of the sample marginals, $\hat{\bm{p}}_{H_0}$. In that case, the substituted term is no longer a constant known value, and convergence should be revisited. There exists generalisations of the $\delta$-method that allow to substitute $\bm{p}_{H_0}$ by a statistic that converges to $\bm{p}_{H_0}$. A sufficient condition is that the derivatives of MI are uniformly differentiable \cite[Theorem 3.9.5]{vaart1996}. Since MI has continuous derivatives in a neighbourhood of $\bm{p}_{H_0}$, it is uniformly differentiable. See also \cite[Note on pg. 386 and Eq. 6a.2.5]{rao1973} for generalisations of the $\delta$-method to non-constant but convergent parameters, and \cite[Note 14.2. on pg. 594]{agresti2012} for generalisation of the $\delta$-method to higher-orders when first derivatives vanishes.
%
%
Notice that marginals convergence, $\hat{\bm{p}}_{H_0} \to \bm{p}_{H_0}$, is faster than joint convergence $\hat{\bm{p}}_n \to \bm{p}$ \cite{agresti2012}.
%
%
%
Note also that $T_n^2$ requires more information than $T_n^1$ in the sense that it uses both the joint distribution and the marginal one in its definition.



Finally, notice that we have used three approximations to derive the test: 1) the effect of the remainder $R(|| \hat{\bm{p}}_n-  \bm{p}_{H_0}||^2)$ has been neglected, 2) considering normality of $\sqrt{n}(\hat{\bm{p}}_n- \bm{p})$ which is only asymptotic, and 3) approximating $\bm{p}_{H_0}$ by $\hat{\bm{p}}_{H_0}$ for $T_n^2$ .



\subsection{Connections with other independence tests}

Detecting statistical dependence has been widely treated in the literature. One of the seminal works on the topic relies on the work of Pearson in 1900, who described the chi-square goodness of fit test based on a geometric approach \cite{pearson1900}. 
Another well-known statistic is the likelihood ratio test presented in 1933 by Neyman and Pearson's son \cite{neyman1933, cressie1989}.
These two are the best-known statistics for testing independence, and their presence is standard in categorical data analysis books.

The likelihood ratio test is based on the ratio \cite[Section 3.2.]{agresti2012}:
\begin{equation}
\Lambda=\frac{\prod_i \prod_j\left(n_{i*} n_{*j}\right)^{n_{i j}}}{n^n \prod_i \prod_j n_{i j}^{n_{i j}}}_. 
\end{equation}

This computes the ratio between the multinomial likelihoods of the hypothesised counts for each cell under $H_0$ and the multinomial likelihood of the actual counts. It is known that the asymptotic distribution of $G^2=-2\ln{\Lambda}$ is $\chi^2_{(I-1)(J-1)}$ \cite[Section 14.3.4.]{agresti2012}. Rearranging the previous equation, we get the typical description of $G_n^2$,
\begin{equation}
G_n^2=-2 \ln \Lambda=2 \sum_i \sum_j n_{i j} \log \left(n_{i j} / ((n_{i*} n_{*j})/n) \right).
\end{equation}
%
Dividing by $n$ on the numerator and denominator inside the logarithm, and dividing and multiplying by $n$ outside, we get the expression of $G_n^2$ in terms of the empirical distribution

\begin{equation}
G_n^2= 2n\sum_i \sum_j \hat{p}_{i j} \ln \left(\hat{p}_{i j} / (\hat{p}_{i*} \hat{p}_{*j}) \right) = 2n\MI(\hat{\bm{p}}_n) = T_n^1.
\end{equation}

which is exactly the statistic $T_n^1$ derived from the MI. Thus, performing the test with $T_n^1$ is equivalent to performing an independence log-likelihood ratio test.

On the other hand, a  interesting connection appears between $T_n^2$ and the classical chi-square test for independence,
\begin{equation}
    \chi_n^2=\sum_{ij=1}^{IJ} \frac{\left(n\hat{p}_{ij}-n \hat{p}_{i*}\hat{p}_{*j}\right)^2}{n \hat{p}_{i*}\hat{p}_{*j}}_.
\end{equation}


We show in Appendix~\ref{ap:godCalculus} that for the plug-in version of $T_{n}^2 $, its expression reduces to $ \chi_n^2$. 
This implies that,

\begin{equation}\label{eq:explicit}
    G^2_n = \chi_n^2 + 2nR(|| \hat{\bm{p}}_n-  \bm{p}_{H_0}||^2)
\end{equation}

where the remainder term can be explicitly quantified by computing third order derivatives of MI. Therefore, equation~\ref{eq:explicit} provides an explicit higher-order relationship between the $G^2_n$ and $\chi_n^2$ statistics.


%

Ultimately,  MI, through its fundamental property that $\MI=0$ if and only if $X$ and $Y$ are independent, establishes an elegant connection between two classical independence test statistics. This finding further emphasises the central role of MI in understanding and quantifying statistical dependence.

\newpage

\bibliographystyle{IEEEtran}
\bibliography{bib}

\newpage


\appendix

\section{Demonstration that $T_n^2$ is the Pearson's chi-square statistic}\label{ap:godCalculus}

Following the notation along this paper, we show that 

$$  T_n^2= \sqrt{n}(\hat{\bm{p}}_n-  \bm{p})^\top H \sqrt{n}(\hat{\bm{p}}_n-  \bm{p})   = \chi_n^2 $$

We split the demonstration by expanding the quadratics terms of each summand of $H = A + B +C$.

\subsection{Quadratic expansion of A}

We compute
$$( \hat{\bm{p}}_n - \hat{\bm{p}}^0_n )^\top A ( \hat{\bm{p}}_n - \hat{\bm{p}}^0_n ). $$

For the sake of simplicity, we remove the hat notation for scalars to alleviate notation. The following calculus stands for any valid probability distribution $\hat{\bm{p}}_n$ with marginals $\hat{\bm{p}}^0_n$.

%

Suppose $i=1, j=1$. From eq. \ref{eq:H}
\begin{equation*}
    A_{1,:} = \begin{pmatrix}
        \frac{1}{p_{ij}}-\frac{1}{p_{i*}} -\frac{1}{p_{*j}}, & -\frac{1}{p_{*j}} & \cdots -\frac{1}{p_{*j}}, & -\frac{1}{p_{i*}} & 0 & \cdots & 0 & 0, & \cdots, & -\frac{1}{p_{i*}} & 0 & \cdots & 0
    \end{pmatrix}_{1 \times IJ-1}
\end{equation*}

Then, 
\begin{align*}
    A_{1,:} \cdot \hat{\bm{p}}^0_n &=  (1 - \frac{1}{p_{1*}} - \frac{1}{p_{*1}}) -  \frac{1}{p_{*1}}\sum_{k=2}^I p^{0}_{k1} -  \frac{1}{p_{1*}}\sum_{k = 2}^J p^{0}_{1k} \\
    & = (1 - \frac{1}{p_{1*}} - \frac{1}{p_{*1}}) - \frac{p_{*1} -p_{11}^0}{p_{*1}}  - \frac{p_{1*} -p_{11}^0}{p_{*1}} \\
    &= -1 - p_{1*} - p_{*1} + \frac{p_{11}^0}{p_{1*}} + \frac{p_{11}^0}{p_{*1}} = -1
\end{align*}

The rows of A corresponding to $i=I$ or $j=J$ has a slightly different value. It can be seen that in general,
\begin{equation*}
    A_{ij,:} \cdot \hat{\bm{p}}^0_n = -1 + \begin{cases}
        0  & i \neq I, \ j \neq  J \\
        \frac{p^0_{IJ}}{p_{I*}} = p_{*J}   & i = I, \ j \neq  J \\
        \frac{p^0_{IJ}}{p_{^J}} = p_{I*}   & i \neq I, \ j =  J
    \end{cases}
\end{equation*}

Operating in a similar manner,

\begin{equation*}
     A_{ij,:} \cdot \hat{\bm{p}}_n = \frac{p_{ij} -2p_{ij}^0}{p_{ij}^0}  + \begin{cases}
        0  & i \neq I, \ j \neq  J \\
        \frac{p_{IJ}}{p_{I*}}   & i = I, \ j \neq  J \\
        \frac{p_{IJ}}{p_{^J}}    & i \neq I, \ j =  J
    \end{cases}
\end{equation*}

Now we multiply by $\hat{\bm{p}}^0_n$ or $\hat{\bm{p}}_n$ on the left side, to be able to evaluate the quadratic form 
$$( \hat{\bm{p}}_n - \hat{\bm{p}}^0_n )^\top A ( \hat{\bm{p}}_n - \hat{\bm{p}}^0_n ) =  \hat{\bm{p}}_n^\top A \hat{\bm{p}}_n -2 (\hat{\bm{p}}_n^0)^\top A   \hat{\bm{p}}_n + (\hat{\bm{p}}_n^0)^\top A \hat{\bm{p}}^0_n. $$

Notice symmetry of A is used. 
With careful consideration we arrive at:
\begin{itemize}
    \item $\hat{\bm{p}}_n^\top A \hat{\bm{p}}_n = \displaystyle\sum_{\substack{i,j=1 \\ i \neq I, j \neq J}}^{I,J}\frac{p_{ij}^2 -2p_{ij}p_{ij}^0}{p_{ij}^0} + (2p_{IJ}-\frac{p_{IJ}^2}{p_{I*}} - \frac{p_{IJ}^2}{p_{*J}})$
    \item $-2 (\hat{\bm{p}}_n^0)^\top A   \hat{\bm{p}}_n = -2 \displaystyle\sum_{\substack{i,j=1 \\ i \neq I, j \neq J}}^{I,J} (p_{ij} - 2p_{ij}^0) -2p_{IJ}(2-p_{I*}-p_{*J})$
    \item $(\hat{\bm{p}}_n^0)^\top A \hat{\bm{p}}^0_n = -\displaystyle\sum_{\substack{i,j=1 \\ i \neq I, j \neq J}}^{I,J}p_{ij}^0  + [2p_{IJ}^0 -(p_{I*} + p_{*J})p_{IJ}^0]$.
\end{itemize}

Notice that with the terms in the summatory we can easily operate to construct the Pearson's chi-square statistic up to term $(I-1,J)$. Putting them together, and simplifying we get:

\[
( \hat{\bm{p}}_n - \hat{\bm{p}}^0_n )^\top A ( \hat{\bm{p}}_n - \hat{\bm{p}}^0_n ) = \sum_{\substack{i,j=1 \\ i \neq I, j \neq J}}^{I,J} \frac{(p_{ij}-p_{ij}^0)^2}{p_{ij}^0}  + \Big[(p_{I*} + p_{*J})(2p_{IJ}-p_{IJ}^0) -\frac{p_{IJ}^2}{p_{I*}} - \frac{p_{IJ}^2}{p_{*J}}\Big]
\]

\subsection{Quadratic expansion of B}

We seek to compute now: 
$$( \hat{\bm{p}}_n - \hat{\bm{p}}^0_n )^\top B ( \hat{\bm{p}}_n - \hat{\bm{p}}^0_n ) $$

The matrix B has identically rows. We exploit this property and multiply by the left to obtain identical terms:

\[
(\hat{\bm{p}}_n^0)^\top B =  \frac{\sum_{k=1}^{J-1} p_{Ik}^0}{p_{I*}}  + \frac{\sum_{k=1}^{I-1} p_{kJ}^0}{p_{I*}} = \frac{p_{I*}-p_{IJ}^0}{p_{I*}}  + \frac{p_{*J}-p_{IJ}^0}{p_{*J}} = 2-p_{*J} -p_{I*}, \ \ \forall \text{ columns }
\]

Analogously,
\[
(\hat{\bm{p}}_n)^\top B = 2-  \frac{ p_{IJ} }{p_{I*}} - \frac{ p_{IJ}}{p_{*J}}, \ \ \forall \text{ columns }
\]

Since the previous vectors have constant terms the quadratics are just the sum of the probabilities $p_{ij}$, except $p_{IJ}$, multiplied by the corresponding constant term:
\begin{itemize}
    \item $(\hat{\bm{p}}_n^\top B ) \hat{\bm{p}}_n  = (2-  \frac{ p_{IJ} }{p_{I*}} - \frac{ p_{IJ}}{p_{*J}})\displaystyle\sum_{\substack{i,j=1 \\ i \neq I, j \neq J}}^{I,J}  p_{ij}  = (2-  \frac{ p_{IJ} }{p_{I*}} - \frac{ p_{IJ}}{p_{*J}})(1-p_{IJ})  $
    \item $((\hat{\bm{p}}_n^0)^\top B) \hat{\bm{p}}^0_n = (2-p_{*J} -p_{I*})(1-p_{IJ}^0)$
    \item $-(\hat{\bm{p}}_n^\top B)   \hat{\bm{p}}_n^0 = - (2-  \frac{ p_{IJ} }{p_{I*}} - \frac{ p_{IJ}}{p_{*J}}) (1-p_{IJ}^0)$
    \item $-((\hat{\bm{p}}_n^0)^\top B)   \hat{\bm{p}}_n = -(2-p_{*J} -p_{I*})(1-p_{IJ})$
\end{itemize}

Notice B is not symmetric.

\subsection{Quadratic expansion of C}

Finally, we seek to compute: 
$$( \hat{\bm{p}}_n - \hat{\bm{p}}^0_n )^\top C ( \hat{\bm{p}}_n - \hat{\bm{p}}^0_n ) $$

The matrix C is identical by columns. We exploit this property and multiplying by the right to obtain identical terms:

\begin{align*}
    C\hat{\bm{p}}_n^0 &= \sum_{l=1}^{J-1} \alpha^0 \sum_{k=1}^{I-1}p_{kl}^0  + \beta^0 \sum_{l=1}^{J-1} p_{Il}  + \gamma^0 \sum_{k=1}^{I-1} p_{kJ}  \\ 
    & = \underbrace{\sum_{l=1}^{J-1} \alpha^0 (p_{*l}-p_{Il}^0)}_{\alpha^0(1-p_{*J}-p_{I*} + p_{IJ}^0)}  +  \beta^0 (p_{I*}-p_{IJ}^0)  +  \gamma^0(p_{*J}-p_{IJ}^0) \ \ \forall \text{ rows }
\end{align*}

Analogously,

\[
 C\hat{\bm{p}}_n =  \alpha^0(1-p_{*J}-p_{I*} + p_{IJ}) +  \beta^0 (p_{I*}-p_{IJ})  +  \gamma^0(p_{*J}-p_{IJ}) \ \ \forall \text{ rows }
\]

This gives four cross product terms of the quadratic form which are:

\begin{itemize}
    \item $\hat{\bm{p}}_n^\top (C  \hat{\bm{p}}_n)  =(1-p_{IJ}) C\hat{\bm{p}}_n$

    \item $(\hat{\bm{p}}_n^0)^\top  (C  \hat{\bm{p}}_n^0) = (1-p_{IJ}^0) C\hat{\bm{p}}_n^0$

    \item $-\hat{\bm{p}}_n^\top (C   \hat{\bm{p}}_n^0) = (1-p_{IJ}) C\hat{\bm{p}}_n^0$

    \item $-(\hat{\bm{p}}_n^0)^\top (C   \hat{\bm{p}}_n) = (1-p_{IJ}^0)C\hat{\bm{p}}_n$
    
\end{itemize}

Notice C is not symmetric.

\subsection{Writing together the terms with B and C}

Now we write the terms with $(1-p_{IJ}^0)$ of the quadratic form

$$( \hat{\bm{p}}_n - \hat{\bm{p}}^0_n )^\top (A + B) ( \hat{\bm{p}}_n - \hat{\bm{p}}^0_n ) $$

together, getting

\begin{align*}
    (1-p_{IJ}^0)\big[ \cancel{2}-p_{*J} -p_{I*} - (\cancel{2}-  \frac{ p_{IJ} }{p_{I*}} - \frac{ p_{IJ}}{p_{*J}})  \big] &+ \alpha^0(1- p_{I*} - p_{*J} + p_{IJ}^0) + \beta^0 (p_{I*} - p_{IJ}^0) + \gamma^0 (p_{*J} - p_{IJ}^0) \\
    &- \alpha^0(1- p_{I*} - p_{*J} + p_{IJ}^0) - \beta^0 (p_{I*} - p_{IJ}^0) - \gamma^0 (p_{*J} - p_{IJ}^0) 
\end{align*}

After careful simplification and direct algebraic manipulation, we arrive at:

\[
\frac{(1-p_{IJ}^0)(1 +p_{*J} + p_{I*})(p_{IJ}-p_{IJ}^0)}{p_{IJ}^0}
\]
By symmetry, the other terms of the quadratic form with  $(1-p_{IJ})$ factor result in:

\[
-\frac{(1-p_{IJ})(1 +p_{*J} + p_{I*})(p_{IJ}-p_{IJ}^0)}{p_{IJ}^0}
\]

Thus, writing them together and taking common factors, 

\[
( \hat{\bm{p}}_n - \hat{\bm{p}}^0_n )^\top (A + B) ( \hat{\bm{p}}_n - \hat{\bm{p}}^0_n ) = \frac{(1+p_{I*}+p_{*J})(p_{IJ}-p_{IJ}^0)^2}{p_{IJ}^0}
\]

\subsection{Final result}

Putting everything together, we get:

\begin{align*}
( \hat{\bm{p}}_n - \hat{\bm{p}}^0_n )^\top (A + B + C) ( \hat{\bm{p}}_n - \hat{\bm{p}}^0_n ) = \sum_{\substack{i,j=1 \\ i \neq I, j \neq J}}^{I,J} \frac{(p_{ij}-p_{ij}^0)^2}{p_{ij}^0}  & + \Big[\cancel{(p_{I*} + p_{*J})(2p_{IJ}-p_{IJ}^0) -\frac{p_{IJ}^2}{p_{I*}} - \frac{p_{IJ}^2}{p_{*J}}}\Big] \\
& +\frac{(1 + \cancel{p_{I*}+p_{*J}})\cancel{(p_{IJ}-p_{IJ}^0)^2}}{p_{IJ}^0}
\end{align*}

For the last cancellation, notice that $-\frac{p_{IJ}^2}{p_{I*}} - \frac{p_{IJ}^2}{p_{*J}} = \frac{-p_{IJ}^2p_{*J}-p_{IJ}^2p_{I*}}{p_{IJ}^0} = \frac{-p_{IJ}^2(p_{I*} + p_{*J})}{p_{IJ}^0}$, thus

\begin{align*}
      \Big[ (p_{I*} + p_{*J})(2p_{IJ}-p_{IJ}^0) -\frac{p_{IJ}^2}{p_{I*}} - \frac{p_{IJ}^2}{p_{*J}} \Big]&=  \frac{(p_{I*} + p_{*J})[2p_{IJ}\cdot p_{IJ}^0 - (p_{IJ}^0)^2 - p_{IJ}^2]}{p_{IJ}^0} \\
      &= -\frac{(p_{I*} + p_{*J}) (p_{IJ} - p_{IJ}^0)^2 }{p_{IJ}^0}
\end{align*}

Then finally,

\[ 
( \hat{\bm{p}}_n - \hat{\bm{p}}^0_n )^\top (A + B + C) ( \hat{\bm{p}}_n - \hat{\bm{p}}^0_n ) = \sum_{\substack{i,j=1 }}^{I,J} \frac{(p_{ij}-p_{ij}^0)^2}{p_{ij}^0}  = \chi_n^2 \qed 
\]

\end{document}